\documentclass[12pt, reqno]{amsart}
\usepackage[utf8]{inputenc}
\usepackage{amsmath,amssymb,amsthm, mathrsfs}
\usepackage{indentfirst}
\usepackage{comment,relsize,braket}
\usepackage{enumerate}
\usepackage{graphicx}
\usepackage{color}
\usepackage[all]{xy}
\usepackage{fullpage}
\usepackage{tikz}
\usepackage[colorlinks, linkcolor=blue!70!black]{hyperref}
\usepackage[top=1in, bottom=1in, left=1in, right=1in]{geometry}
\usepackage{times}
\usepackage{tfrupee}
\usetikzlibrary{chains,fit}
\usetikzlibrary{shapes,snakes}
\usetikzlibrary{graphs}
\usepackage{graphicx,adjustbox}
\usepackage{hyperref}
\usepackage{amscd}
\usepackage{booktabs}

\usepackage{float}
\usepackage{caption}
\usepackage{subcaption}

\usepackage{cleveref}

\newcommand{\cH}{\mathcal{H}}

\theoremstyle{plain}

\theoremstyle{theorem}

\newtheorem{definition}{Definition}[section]
\newtheorem{proposition}[definition]{Proposition}
\newtheorem{theorem}[definition]{Theorem}
\newtheorem{lemma}[definition]{Lemma}
\newtheorem{conjecture}[definition]{Conjecture}

\newtheorem{example}[definition]{Example}

\newtheorem{remark}[definition]{Remark}

\theoremstyle{remark}

\makeatletter
\@namedef{subjclassname@2020}{\textup{2020} Mathematics Subject Classification}
\makeatother

\begin{document}
		
		\title[Symmetric generalized numerical semigroups in $\mathbb{N}^d$ with embedding dimension $2d+1$]{Symmetric generalized numerical semigroups in $\mathbb{N}^d$ with embedding dimension $2d+1$}

         \author{OM PRAKASH BHARDWAJ}
		\address{Chennai Mathematical Institute, Siruseri, Chennai, Tamil Nadu 603103, India
			}
		\email{omprakash@cmi.ac.in; om.prakash@alumni.iitgn.ac.in}        
  
		\author{CARMELO CISTO}
		\address{Universit\`{a} di Messina, Dipartimento di Scienze Matematiche e Informatiche, Scienze Fisiche e Scienze della Terra\\
			Viale Ferdinando Stagno D'Alcontres 31\\
			98166 Messina, Italy}
		\email{carmelo.cisto@unime.it}

		

		\begin{abstract}

   In this article, we classify all symmetric generalized numerical semigroups in $\mathbb{N}^d$ of embedding dimension $2d+1$.   Consequently, we show that in this case the property of being symmetric is equivalent to have a unique maximal gap with respect to natural partial order in $\mathbb{N}^d$. Moreover, we deduce that when $d>1$, there does not exist any generalized numerical semigroup of embedding dimension $2d+1$, which is almost symmetric but not symmetric.
  	% In this article, we classify all symmetric generalized numerical semigroups of embedding dimension $2d+1$.
   % %\textcolor{blue}{Consequently, when $d>2$, if a generalized numerical semigroup is symmetric, then its embedding dimension can not be $2d+1$}. 
   % Consequently, we deduce that when $d>1$, there does not exist any generalized numerical semigroup of embedding dimension $2d+1$, which is almost symmetric but not symmetric.

		\end{abstract}

		\maketitle

\section*{Introduction}

Let $\mathbb{N}$ be the set of nonnegative integers. A submonoid of $\mathbb{N}$ having finite complement in it is called a \emph{numerical semigroup}. This notion is widely studied in the scientific literature, starting from its first appearance in \cite{sylvester1884mathematical}, where the so called \emph{Frobenius problem} is introduced. The research in this topic is still active and one can refer to \cite{rosales2009numerical} and \cite{numericalApplications} for good expository references about the main properties of these monoids, together with a collection of many results obtained in the development of this line of research. Numerical semigroups are also related with other wide topics, like Commutative Algebra and Algebraic Geometry. 

%\textcolor{red}{ About commutative algebra.... For symmetric related to Gorenstein: \cite{kunz1970value}. For almost-symmetric and almost Gorenstein: \cite{barucci1997one}. The concept of type}

The concept of numerical semigroups can be generalized considering submonoids of $\mathbb{N}^d$ with finite complement in it, for $d$ any positive integer. This kind of generalization is considered in \cite{failla2016algorithms}, introducing the notion of \emph{generalized numerical semigroup}. Some properties of these monoids have been introduced in some recent papers. For example, it is known that every generalized numerical semigroup is finitely generated and admits a unique minimal set of generators that is finite (see \cite{cisto-generators}). So, if $S$ is a generalized numerical semigroup in $\mathbb{N}^d$, $d\geq 1$, it is possible to define the invariant $\operatorname{e}(S)$, consisting of the cardinality of its minimal set of generators, called \emph{embedding dimension} of $S$. Considering the wide scientific literature about numerical semigroups, a natural direction of research is to study how some properties and notions can be extended in the context of generalized numerical semigroups. In this direction, the notions of symmetric and almost symmetric generalized numerical semigroup are introduced and studied, respectively, in \cite{cisto-irreducible} and \cite{cisto-tenorio} (see also \cite{op1,garcia-pseudofrobenius} for an extension of this notions in an even more general framework).

It is well known that every numerical semigroup of the embedding dimension $2$ is symmetric (see \cite[Proposition 2.13 and Corollary 4.5]{rosales2009numerical}). Moreover, for numerical semigroups of embedding dimension $3$, a classification of almost symmetric property is provided by \cite{Herzog} (for symmetric) and by \cite{pedroEmbedding3} (for not symmetric). Since for a generalized numerical semigroup properly contained in $\mathbb{N}^d$, the least possible value of embedding dimension is $2d$, a natural question is to classify the symmetric and almost symmetric generalized numerical semigroups with embedding dimension $2d$ and $2d+1$. In the case of embedding dimension equal to $2d$, such a classification can be recovered by some results in \cite{cisto2020generalization}. The investigation of the case of embedding dimension $2d+1$ has not been carried out till now and this is the purpose of this paper. 

The paper is structured as follows. In Section 1, we recall the main notations and concepts related to generalized numerical semigroups, that will be useful throughout the paper. Moreover, we state in Proposition~\ref{prop:2d} the result for the classification of symmetric generalized numerical semigroups in $\mathbb{N}^d$ with embedding dimension $2d$. In Section 2, we provide some preliminary results that will allow to obtain the main result of this paper. In particular, in Lemma~\ref{lem:class1}, \ref{lem:class2} and \ref{lem:class3}, we classify the generalized numerical semigroups in $\mathbb{N}^d$ of embedding dimension $2d+1$ having unique Frobenius allowable element. From these results, in Section 3 we prove Theorem~\ref{thm:symmteric2d+1}, which states that a generalized numerical semigroup in $\mathbb{N}^d$ of embedding dimension $2d+1$ is symmetric if and only if it has a unique Frobenius allowable element. Consequently, when $d>1$, we argue that symmetric and almost symmetric properties are equivalent in embedding dimension $2d+1$. In particular, for every positive integer $d$, we can state that for an almost symmetric generalized numerical semigroup in $\mathbb{N}^d$ of embedding dimension $2d+1$, the type is bounded by $2$. Finally, encouraged by some computational experiments, we provide a conjecture about the type of an almost symmetric generalized numerical semigroup of embedding dimension $2d+2$ (see Conjecture~\ref{conj:2d+2}).

\section{Some notions on generalized numerical semigroups}

Let $\mathbb{N}$ be the set of nonnegative integers, and $d$ a positive integer. In the set $\mathbb{N}^d$ we define $\leq_{\mathbb{N}^d}$ to be the natural partial order, that is, for $\mathbf{x},\mathbf{y}\in \mathbb{N}^d$ we have $\mathbf{x}\leq_{\mathbb{N}^d} \mathbf{y}$ if $\mathbf{y}-\mathbf{x} \in \mathbb{N}^d$. Moreover, we denote by $\mathbf{e}_1,\ldots,\mathbf{e}_d$ the standard basis vectors of the vector space $\mathbb{R}^d$.

An additive submonoid $S$ of $\mathbb{N}^d$ is called a generalized numerical semigroup if the set $\mathcal{H}(S)=\mathbb{N}^d\setminus S$ is finite. In particular, the set $\mathcal{H}(S)$ is called the set of \textit{gaps} of $S$. When $d=1$, $S$ is a numerical semigroup and the largest integer belonging to $\mathcal{H}(S)$ is called the Frobenius number of the numerical semigroup $S$, denoted by $\mathrm{F}(S)$.  If $A\subseteq \mathbb{N}^d$, we denote $\langle A \rangle=\{\sum_{i=1}^n \lambda_i \mathbf{a}_i\mid n \in \mathbb{N}\setminus \{0\}, \lambda_i \in \mathbb{N}\setminus \{0\}, \mathbf{a}_i \in A, \text{ for all } i\in \{1,\ldots,n\}\}$, that is the submonoid of $\mathbb{N}^d$ generated by $A$. By \cite{cisto-generators}, it is known that for every generalized numerical semigroup $S$ there exists a unique finite minimal set of generators of it, that is, there exists a unique finite set $A$ such that $S=\langle A \rangle$ and for all $B\subsetneq A$ we have $S\neq \langle B \rangle$. The cardinality of the minimal set of generators is called \textit{embedding dimension} and denoted by $\operatorname{e}(S)$. As shown in \cite[Theorem 11]{garcia2018extension}, the value of the embedding dimension of a generalized numerical semigroup has the following lower bound:

\begin{proposition}\cite{garcia2018extension}
Let $S\subsetneq \mathbb{N}^d$ be a generalized numerical semigroup. Then $\operatorname{e}(S)\geq 2d$.
\end{proposition}

If $S\subseteq \mathbb{N}^d$ is a generalized numerical semigroup, the following sets are useful for the purpose of this paper:

\begin{itemize}
\item $\operatorname{PF}(S)=\{\mathbf{h}\in \mathcal{H}(S)\mid \mathbf{h}+\mathbf{s} \in S,\text{ for all }\mathbf{s}\in S\setminus \{\mathbf{0}\}\}$ is called the set of \textit{pseudo-Frobenius} elements. We also denote $\operatorname{t}(S)=|\operatorname{PF}(S)|$, called the \textit{type} of $S$.

\item $\operatorname{FA}(S)=\mathrm{Maximals}_{\leq_{\mathbb{N}^d}} \mathcal{H}(S)$ is called the set \textit{Frobenius-allowable} elements (see also \cite{singhal} for the meaning of this notion).
\end{itemize}

\begin{example} \rm
Let $S = \mathbb{N}^2\setminus  \{(0, 1), (1, 0), (1, 1), (1, 2), (1, 3), (1, 4), (2, 1), (3, 0), (3, 2)\}$. It is not difficult to see that $S$ is a generalized numerical semigroup with $\operatorname{FA}(S)=\{(1,4),(3,2)\}$, and $\operatorname{PF}(S)=\{(1,4),(3,2),(1,3),(2,1)\}.$
\end{example}

It is easy to see that $\operatorname{FA}(S)\subseteq \operatorname{PF}(S)$. These sets allow us to characterize the notions of symmetric and almost symmetric generalized numerical semigroup. These concepts have also been introduced and studied for more general classes of finitely generated monoids in $\mathbb{N}^d$ (see, for instance, \cite{op1, garcia-pseudofrobenius, garcia_someproperties}). For simplicity, we provide these notions in the case of generalized numerical semigroups. 

If $S\subseteq \mathbb{N}^d$ is a generalized numerical semigroup, for our purpose we can consider the following equivalent formulations of symmetric and almost symmetric, that one can find in \cite[Definition 2.11]{cisto-irreducible} and \cite[Propositions 2.14 and 3.3]{cisto-tenorio}:

%\begin{definition}\rm
\begin{enumerate}
\item $S$ is symmetric if and only if $\operatorname{t}(S)=1$.

\item $S$ is almost symmetric if and only if $\operatorname{FA}(S)=\{\mathbf{f}\}$ such that $\mathbf{f}-\mathbf{p} \in \operatorname{PF}(S)$ for all $\mathbf{p}\in \operatorname{PF}(S)\setminus \{\mathbf{f}\}$.
\end{enumerate}
In particular, we can say that a symmetric generalized numerical semigroup is an almost symmetric generalized numerical semigroup such that $|\operatorname{PF}(S)|=1$. Furthermore, an almost symmetric generalized numerical semigroup having $\operatorname{t}(S)=2$ is usually called \textit{pseudo-symmetric}. 
%\end{definition} 

%\subsection*{Commutative algebraic correspondence:}

\medskip
Now, we recall a result which gives a criterion for the symmetric property of a generalized numerical semigroup.

\begin{theorem}[{\cite[Theorem 5.6]{cisto-irreducible}}]\label{symcrit}   
Let $S \subseteq \mathbb{N}^d$ be a generalized numerical semigroup with $\vert \mathcal{H}(S) \vert = g$. Then $S$ is symmetric if and only if there exists $\mathbf{f} = (f_1, \ldots, f_d) \in \mathcal{H}(S)$ with $2g = ( f_1 + 1)(f_2+1) \cdots ( f_d + 1)$. 
\end{theorem}

\medskip
Let us consider the set of generalized numerical semigroups in $\mathbb{N}^d$ with embedding dimension $2d$ (that is, the least possible value). If $d=1$, then it is well known that every numerical semigroup generated by two integers is symmetric (see, for instance \cite[Proposition 2.13 and Corollary 4.5]{rosales2009numerical}). In the case $d>1$, by the proof of Lemma 4.1 and Remark 4.2 of \cite{cisto2020generalization}, we have the following result.

\begin{proposition}\label {prop:2d}
Let $d>1$ and $S\subseteq \mathbb{N}^d$ be a generalized numerical semigroup such that $\operatorname{e}(S)=2d$. Then $S$ is symmetric if and only if $|\operatorname{FA}(S)|=1$. Also, in this case, sets of symmetric and almost symmetric generalized numerical semigroups coincide.

\noindent Moreover, in this case, we have $S=\langle A \rangle$, with $A=\{\mathbf{e}_1, \ldots, \mathbf{e}_{i-1}, 2\mathbf{e}_i, b\mathbf{e}_i, \mathbf{e}_{i+1},\ldots,\mathbf{e}_d\}\cup \{\mathbf{e}_i+d_j\mathbf{e}_j\mid j\in \{1,\ldots,d\}\setminus \{i\}\}\subseteq \mathbb{N}^d$, where $b$ is an odd number and $d_j\in \mathbb{N}\setminus \{0\}$ for all $j\in \{1,\ldots,d\}\setminus \{i\}$.

%By the proof of \cite[Lemma 4.1]{cisto2020generalization}
%$$\cH(S)=\left\lbrace h\mathbf{e}_1+\sum_{j\neq 1} \lambda_j\mathbf{e}_j\mid h\in \cH(T), \lambda_j\in \{0,1,\ldots,d_j-1\}, j\in \{2,\ldots,d\}\right \rbrace.$$

\end{proposition}
\begin{proof}
    The first claim and the last claim are consequence of \cite[Remark 4.2]{cisto2020generalization}. The second claim follows since every symmetric generalized numerical semigroup is almost symmetric, and if $S$ is almost symmetric then $|\operatorname{FA}(S)|=1$.
\end{proof}

 In successive sections, we study the almost symmetric generalized numerical semigroups with embedding dimension $2d+1$

\section{Preparatory results}

If $S\subseteq \mathbb{N}^d$ is a generalized numerical semigroup such that $\operatorname{e}(S)=2d+1$ and $d>1$, by \cite[Theorem 2.8]{cisto-generators}, we have the following possibilities (up to permutation of coordinates):

\begin{enumerate}
    \item $S=\langle A \rangle$, with $A:=\{\mathbf{e}_1, \ldots, \mathbf{e}_{i-1}, a_1\mathbf{e}_i, a_2\mathbf{e}_i, a_3\mathbf{e}_i, \mathbf{e}_{i+1},\ldots,\mathbf{e}_d\}\cup \{\mathbf{e}_i+d_j\mathbf{e}_j\mid j\in \{1,\ldots,d\}\setminus \{i\}\}\subseteq \mathbb{N}^d$, $d>1$ an integer, such that $a_1,a_2,a_3$ minimally generate a numerical semigroup and $d_j\in \mathbb{N}\setminus \{0\}$ for all $j\in \{1,\ldots,d\}\setminus \{i\}$.
    
    %$S_1=\langle (a,0),(b,0),(c,0),(0,1), (1,d) \rangle$ with $d\geq 1$ and $\gcd(a,b,c)=1$.

    \medskip
    \item $S=\langle (a_1,0),(a_2,0),(0,b_1),(0,b_2),(1,1)\rangle\subseteq \mathbb{N}^2$, with $\gcd(a_1,a_2)=1$ and $\gcd(b_1,b_2)=1$.

 \medskip
    \item $S=\langle A \rangle$, with $A:=\{\mathbf{e}_1, \ldots, \mathbf{e}_{i-1}, a_1\mathbf{e}_i, a_2\mathbf{e}_i,\mathbf{e}_{i+1},\ldots,\mathbf{e}_d\}\cup \{\mathbf{e}_i+d_j\mathbf{e}_j\mid  j\in \{1,\ldots,d\}\setminus \{i\}\}\cup \{\mathbf{x}\}\subseteq \mathbb{N}^d$ such that $\gcd(a_1,a_2)=1$, $d_j\in \mathbb{N}\setminus \{0\}$ for all $j\in \{1,\ldots,d\}\setminus \{i\}$, $\mathbf{x}\notin \langle A\setminus \{\mathbf{x}\}\rangle$ and $\mathbf{x}\notin \{\lambda \mathbf{e}_i, \mathbf{e}_i+\beta_j\mathbf{e}_j \mid \lambda,\beta \in \mathbb{N}\}$.
    
    %$S_3=\langle (a,0),(b,0),(0,1), (1,d), \mathbf{x} \rangle$ with $\mathbf{x}=(x_1,x_2) \notin \langle (a,0),(b,0),(0,1), (1,d)\rangle$, $x_1>1$ and $x_2>0$.
\end{enumerate}

% \begin{lemma}
%     Let $S_1=\langle (a,0),(b,0),(c,0),(0,1), (1,d) \rangle$ with $d\geq 1$ and $\gcd(a,b,c)=1$. Then $\vert \operatorname{FA}(S_1)\vert=1$ if and only if $a=3$, $b\equiv 1 \mod 3$ and $c=b+1$.
% \end{lemma}
\begin{lemma}\label{lem:class1}
    Let $d>1$ be an integer and $A:=\{\mathbf{e}_1, \ldots, \mathbf{e}_{i-1}, a_1\mathbf{e}_i, a_2\mathbf{e}_i, a_3\mathbf{e}_i, \mathbf{e}_{i+1},\ldots,\mathbf{e}_d\}\cup \{\mathbf{e}_i+d_j\mathbf{e}_j\mid j\in \{1,\ldots,d\}\setminus \{i\}\}\subseteq \mathbb{N}^d$ such that $a_1,a_2,a_3$ minimally generate a numerical semigroup and $d_j\in \mathbb{N}\setminus \{0\}$ for all $j\in \{1,\ldots,d\}\setminus \{i\}$. Let $S=\langle A \rangle $, then the following holds:
    \begin{enumerate}
         \item If $d\geq 3$ then $|\operatorname{FA}(S)|\geq 2$.
        \item If $d=2$ and $a_1\neq3$ then $|\operatorname{FA}(S)|\geq 2$.
        \item If $d=2$ and $a_1=3$, then 
             $|\operatorname{FA}(S)|=1$ if and only if $a_2\equiv 1 \mod 3$ and $a_3=a_2+1$.
        
    \end{enumerate}
\end{lemma}
\begin{proof}
Without loss of generality, we can assume that $i=1$. Let us denote that $T=\langle a_1,a_2, a_3\rangle$. Since $T$ has three minimal generators, we have $a_1\geq 3$. Suppose first $\operatorname{F}(T)-1\in T$. We show that $\mathbf{f}=\operatorname{F}(T)\mathbf{e}_1+\sum_{j\neq1}(d_i-1)\mathbf{e}_i$ is maximal in $\cH(S)$ with respect to $\leq_{\mathbb{N}^d}$. Assuming $\mathbf{f}\in S$, observe that the generators $\mathbf{e}_1+d_j\mathbf{e}_j$ are not involved in the factorizations of $\mathbf{f}$. Hence, since the first coordinate of $\mathbf{f}$ is nonzero, the only possibility is that $\operatorname{F}(T)\in T$. This is a contradiction. Therefore, we get $\mathbf{f}\in \cH(S)$. In order to show that $\mathbf{f}$ is a maximal gap, we take $\mathbf{y}=\sum_{j=1}^dy_j\mathbf{e}_j\in \mathbb{N}^d\setminus \{\mathbf{0}\}$ and we prove that $\mathbf{f}+\mathbf{y}\in S$. We can set $y_1=0$, otherwise it is trivial. Let $k\in \{2,\ldots,d\}$ such that $y_k>0$. So, we get 
$$\mathbf{f}+\mathbf{y}=(\operatorname{F}(T)-1)\mathbf{e}_1+\sum_{j\neq 1,k}(d_j-1+y_j)\mathbf{e}_j+(\mathbf{e}_1+d_k\mathbf{e}_k)+(y_k-1)\mathbf{e}_k\in S.$$
Hence, $\mathbf{f}$ is maximal in $\cH(S)$. Moreover, for all $k\in \{2,\ldots,d\}$, consider the element $\mathbf{h}_k=(a_1-1)\mathbf{e}_1+(a_1-1)(\mathbf{e}_1+d_k\mathbf{e}_k)-\mathbf{e}_k$. It is not difficult to see that $\mathbf{h}_k\in \cH(S)$ and it is not comparable with $\mathbf{f}$ with respect to $\leq_{\mathbb{N}^d}$. In particular, $|\operatorname{FA}(S)|\geq 2$. Now, suppose $\operatorname{F}(T)-1\notin T$. Observe that, in this case, $\operatorname{F}(T),\operatorname{F}(T)-1\in \operatorname{PF}(T)$ and if $\operatorname{F}(T)-2\notin T$, then $\operatorname{F}(T)-2\in \operatorname{PF}(T)$. But, by \cite[Corollary 10.22]{rosales2009numerical}, we have that $\operatorname{t}(T)\leq 2$. Consequently, we have $\operatorname{F}(T)-2\in T$. For all $k\in \{2,\ldots,d\}$, denote $\mathbf{f}_k=\operatorname{F}(T)\mathbf{e}_1+\sum_{j\neq 1,k}(d_j-1)\mathbf{e}_j+(2d_k-1)\mathbf{e}_k$. Assume $\mathbf{f}_k\in S$. Observe that the generators $\mathbf{e}_1+d_j\mathbf{e}_j$, for $j\in \{2,\ldots,d\}\setminus \{k\}$, are not involved in the factorizations of $\mathbf{f}_k$. If $\mathbf{e}_k+d_k\mathbf{e}_k$ is involved in any factorization of $\mathbf{f}_k$, then 
$$\mathbf{f}_k-(\mathbf{e}_1+d_k\mathbf{e}_k)=(\operatorname{F}(T)-1)\mathbf{e}_1+\sum_{j\neq1}(d_j-1)\mathbf{e}_j\in S,$$
this is a contradiction by similar arguments used before and since $\operatorname{F}(T)-1\notin T$. So, also $\mathbf{e}_1+d_k\mathbf{e}_k$ is not involved in the factorizations of $\mathbf{f}_k$. As a consequence, by the assumption $\mathbf{f}_k  \in S$, we have that $\operatorname{F}(T)\in T$, a contradiction. Hence, $\mathbf{f}_k\in \cH(S)$ for all $k\in \{2,\ldots,d\}$. We want to show that, for all $k\in \{2,\ldots,k\}$, $\mathbf{f}_k$ is maximal in $\cH(S)$ with respect to $\leq_{\mathbb{N}^d}$. Let $\mathbf{y}=\sum_{j=1}^dy_j\mathbf{e}_j\in \mathbb{N}^d\setminus \{\mathbf{0}\}$. If $y_1>0$, then trivially $\mathbf{f}_k+\mathbf{y}\in S$. If there exists $\ell\in \{2,\ldots,d\}\setminus \{k\}$ such that $y_\ell>0$, then 
$$\mathbf{f}_k+\mathbf{y}=(\operatorname{F}(T)-2)\mathbf{e}_1+(\mathbf{e}_1+d_j\mathbf{e}_j)+\sum_{j\neq 1,\ell}(d_j-1+y_j)\mathbf{e}_j+(y_\ell-1)\mathbf{e}_\ell+(\mathbf{e}_1+d_k\mathbf{e}_k)\in S,$$ 
observing that in this case $\operatorname{F}(T)-2\in T$. If $y_k>0$, then 
$$\mathbf{f}_k+\mathbf{y}=(\operatorname{F}(T)-2)\mathbf{e}_1+\sum_{j\neq 1,k}(d_j-1+y_j)\mathbf{e}_j+2(\mathbf{e}_1+d_k\mathbf{e}_k)+ (y_k-1)\mathbf{e}_k \in S.$$ 
Hence, $\mathbf{f}_k$ is maximal in $\cH(S)$ for all $k\in \{2,\ldots,d\}$. Now observe that if $d>2$, then for all $k_1,k_2\in \{2,\ldots,d\}$, with $k_1\neq k_2$, we have $\mathbf{f}_{k_1}$ and $\mathbf{f}_{k_2}$ are not comparable with respect to $\leq_{\mathbb{N}^d}$. In particular, $|\operatorname{FA}(S)|\geq 2$. Therefore, suppose $d=2$. For simplicity, we denote $S=\langle (a_1,0),(a_2,0),(a_3,0), (0,1),(1,d_2) \rangle$ and by the previous arguments $\mathbf{f}=(\operatorname{F}(T), 2d_2-1)$ is maximal in $\cH(S)$. Now, if $a_1\geq 4$, then $(3,3d_2-1)\in \cH(S)$ and it is not comparable to $\mathbf{f}$. In particular, $|\operatorname{FA}(S)|\geq 2$. If $a_1=3$, considering that we are assuming $\operatorname{F}(T)-1\notin T$, we have that $a_2\equiv 1 \mod 3$ and $a_3=a_2+1$. In this case, is not difficult to show that $\cH(S)=\{(h,j)\mid h\in \cH(T), h\equiv 1 \mod 3, 0\leq j\leq d_2-1\}\cup
\{(h,j)\mid h\in \cH(T), h\equiv 2 \mod 3, 0\leq j\leq 2d_2-1\}$. In particular, $\operatorname{FA}(S)=\{(\operatorname{F}(T),2d_2-1)\}$. This completes the proof.
\end{proof}

\begin{lemma}\label{lem:class2}
  Let $S=\langle (a_1,0),(a_2,0),(0,b_1),(0,b_2),(1,1)\rangle$, with $\gcd(a_1,a_2)=1$ and $\gcd(b_1,b_2)=1$. Then $\vert \operatorname{FA}(S)\vert \geq 2$ 
\end{lemma}
\begin{proof}
Let us define the numerical semigroups $T_1=\langle a_1, a_2 \rangle$ and $T_2=\langle b_1, b_2\rangle.$ First suppose $a_1\notin T_2$. In this case we show that $\mathbf{f}_1=(a_1-1,a_1-1+\operatorname{F}(T_2))$ is a maximal gap pf $S$ with respect to $\leq_{\mathbb{N}^2}$. Note that $\mathbf{f_1} \notin S$, because if  $\mathbf{f}_1\in S$,
the only possibility for this to happen is $\mathbf{f}_1 = (a_1-1)(1,1) + \mathbf{s} ~\text{for some}~\mathbf{s}\in S.$ But the only choice for $\mathbf{s}$ is $(0,\mathrm{F}(T_2))$, this is a contradiction since $(0,\mathrm{F}(T_2)) \notin S$.
% Define $$k_1 =
%  \max\{i\in \mathbb{N}\mid ia_1<a_2\} ~\text{and}~ k_2 = \max\{i\in \mathbb{N}\mid ib_1<b_2\}.$$ Consider the vectors $\mathbf{f}_1=(k_1 a_1-1,a_1-1+\operatorname{F}(T_2))$ and $\mathbf{f}_2=(b_1-1+\operatorname{F}(T_1), k_2 b_1-1)$.
% We want to show that $\mathbf{f}_1$ and $\mathbf{f}_2$ are maximal gaps with respect to $\leq_{\mathbb{N}^2}$. We first examine $\mathbf{f}_1$. 
% Assume $\mathbf{f}_1\in S$. Since $k_1 a_1 < a_2$ and $k a_1-a_1\notin T_1$, the only possibility for this to happen is $\mathbf{f}_1 = (k_1a_1-1)(1,1) + \mathbf{s} ~\text{for some}~\mathbf{s}\in S.$ But the only choice for $s$ is $(0,\mathrm{F}(T_2))$, this is a contradiction since $(0,\mathrm{F}(T_2)) \notin S$. Hence, $\mathbf{f}_1\in \cH(S)$.
Now we prove that $\mathbf{f}_1+(x,y)\in S$ for all $(x,y)\in \mathbb{N}^2 \setminus \{\mathbf{0}\}$. Write $x=qa_1+r$ for some $q,r\in \mathbb{N}$ and $r<a_1$. Assuming $r>0$ we have:
\begin{align*}
    \mathbf{f}_1+(x,y) & = (a_1-1, a_1-1+\operatorname{F}(T_2))+(x,0)+(0,y) \\
    & = (a_1-1, a_1-1+\operatorname{F}(T_2))+(qa_1,0)+(r,0)+(0,y) \\
    & = q(a_1,0)+(a_1,a_1-1+\operatorname{F}(T_2)-(r-1)+y)+(r-1)(1,1).
\end{align*}
Since $r-1 < a_1-1$, we get  $\mathbf{f}_1+(x,y) \in S.$ Now, suppose $r=0, y>0$. This implies $x = qa_1$ and
\begin{align*}
    \mathbf{f}_1+(x,y) & = (a_1-1, a_1-1+\operatorname{F}(T_2))+(qa_1,0)+(0,y) \\
    & = q(a_1,0)+(a_1-1)(1,1)+(0,\operatorname{F}(T_2)+y).
\end{align*}
Since $y>0$ then $\operatorname{F}(T_2)+y\in T_2$ and hence $\mathbf{f}_1+(x,y)\in S$. Now suppose $r=0, y=0$. Since $(x,y)\neq (0,0)$ we have $q\geq 1$. In this case, we have
\begin{align*}
    \mathbf{f}_1+(x,y) & = (a_1-1, a_1-1+\operatorname{F}(T_2))+(qa_1,0) \\
    & = (q-1)(a_1,0)+(2a_1-1)(1,1)+(0,\operatorname{F}(T_2)-a_1).
\end{align*}
Since $a_1 \notin T_2$ and $T_2$ is symmetric, we get $\operatorname{F}(T_2)-a_1\in T_2$ (see \cite[Proposition 4.4]{rosales2009numerical}). Therefore, $\mathbf{f}_1+(x,y)\in S$.

\noindent Now suppose $a_1\in T_2$. Define $k_1 =
 \max\{i\in \mathbb{N}\mid ia_1<a_2\}$ and  consider the vector $\mathbf{f}'_1=(k_1 a_1-1,a_1-1+\operatorname{F}(T_2))$. We want to show that in this case, $\mathbf{f}_1'$ is a maximal gap of $S$ with respect to $\leq_{\mathbb{N}^2}$. First, we must prove that $\mathbf{f}_1'\notin S$. Assume $\mathbf{f}_1'\in S$. Since $k_1a_1<a_2$, the only choices  for a factorization of $\mathbf{f}_1'$ are
 $$(k_1-i)(a_1,0)+(ia_1-1)(1,1)+(0,s) ~\text{for some}~ s\in T_2 ~\text{with}~  1\leq i\leq k_1.$$

\noindent But in this way we obtain $s=\operatorname{F}(T_2)-ia_1\notin T_2$, since $a_1\in T_2$. This is a contradiction, so $\mathbf{f}_1'\notin S$.
Now we prove that $\mathbf{f}_1'+(x,y)\in S$ for all $(x,y)\in \mathbb{N}^2 \setminus \{\mathbf{0}\}$. Write $x=qa_1+r$ for some $q,r\in \mathbb{N}$ and $r<a_1$. Assuming $r>0$ we have:

\begin{align*}
    \mathbf{f}_1'+(x,y) & = (k_1 a_1-1, a_1-1+\operatorname{F}(T_2))+(qa_1+r,0)+(0,y) \\
    & = q(a_1,0)+(k_1 a_1,a_1-1+\operatorname{F}(T_2)-(r-1)+y)+(r-1)(1,1).
\end{align*}

So, since $a_1>r$, we obtain $\mathbf{f}_1'+(x,y)\in S$. Assume $r=0, q\geq 1$, then we have:

\begin{align*}
    \mathbf{f}_1'+(x,y) & = (k_1 a_1-1, a_1-1+\operatorname{F}(T_2))+(q a_1,0)+(0,y) \\
    & = (a_2,0)+((k_1+1)a_1-a_2-1)(1,1)+(q-1)(a_1,0)+(0,a_2-k_1a_1+\operatorname{F}(T_2)+y).
\end{align*}

By the defintion of $k_1$ and $a_2$ is a minimal generator of $T_1$, we get $(k_1+1)a_1>a_2$, and we obtain $\mathbf{f}_1'+(x,y)\in S$. Now assume $r=0, q=0$, then we have $y>0$ and we obtain:

\begin{align*}
    \mathbf{f}_1'+(x,y) & = (k_1a_1-1, a_1-1+\operatorname{F}(T_2))+(0,y) \\
    & = (a_1-1)(1,1)+(k_1-1)(a_1,0)+(0,\operatorname{F}(T_2)+y).
\end{align*}
Since $y>0$ then $\operatorname{F}(T_2)+y\in T_2$ and hence $\mathbf{f}'_1+(x,y)\in S$. Therefore we proved that, in any case, $\mathbf{f}_1'+(x,y)\in S$ for all $(x,y)\in \mathbb{N}^2\setminus \{(0,0\}$, that is, $\mathbf{f}_1'\in \operatorname{FA}(S)$.

\noindent With the same arguments, we can prove that if $b_1 \notin T_1$, then $\mathbf{f}_2=(b_1-1+\operatorname{F}(T_1),b_1-1)$ is a maximal gap of $S$ with respect to $\leq_{\mathbb{N}^2}$. While, if $b_1\in T_1$, define $k_2 = \max\{i\in \mathbb{N}\mid ib_1<b_2\}$, and the vector $\mathbf{f}_2'=(b_1-1+\operatorname{F}(T_1), k_2 b_1-1)$ is a maximal gap of $S$ with respect to $\leq_{\mathbb{N}^2}$.

\noindent Now, in order to complete the proof, it is sufficient to show that $\mathbf{f}_1\nleq_{\mathbb{N}^2} \mathbf{f}_2$, $\mathbf{f}_1\nleq_{\mathbb{N}^2} \mathbf{f}_2'$, $\mathbf{f}_1'\nleq_{\mathbb{N}^2} \mathbf{f}_2$ and $\mathbf{f}'_1\nleq_{\mathbb{N}^2} \mathbf{f}_2'$. It easily follows that $a_1-1<b_1-1+\operatorname{F}(T_1)$ and $b_1-1<a_1-1+\operatorname{F}(T_2)$. This implies $\mathbf{f}_1\nleq_{\mathbb{N}^2} \mathbf{f}_2$. Moreover, since $k_2b_1-1\notin T_2$ and $a_1\geq 2$, we get $k_2b_1-1<\operatorname{F}(T_2)+a_1-1$. Hence, we get $\mathbf{f}_1\nleq_{\mathbb{N}^2} \mathbf{f}_2'$. With the same arguments, we can show that $k_1a_1-1<b_1-1+\operatorname{F}(T_1)$. So, together with the inequalities obtained above, we also obtain $\mathbf{f}_1'\nleq_{\mathbb{N}^2} \mathbf{f}_2$ and $\mathbf{f}'_1\nleq_{\mathbb{N}^2} \mathbf{f}_2'$. This completes the proof.
% If $r>0$, then $r-1 \geq 0$ and $a_1-1+\operatorname{F}(T_2)-(r-1)+y \in T_2$ because $r-1 < a_1-1$, hence $\mathbf{f}_1+(x,y)\in S$. Now suppose $r=0$, 
\end{proof}

\begin{lemma} \label{lem:class3}
    Let $d>1$ be an integer and $A:=\{\mathbf{e}_1, \ldots, \mathbf{e}_{i-1}, a_1\mathbf{e}_i, a_2\mathbf{e}_i,\mathbf{e}_{i+1},\ldots,\mathbf{e}_d\}\cup \{\mathbf{e}_i+d_j\mathbf{e}_j\mid  j\in \{1,\ldots,d\}\setminus \{i\}\}\cup \{\mathbf{x}\}\subseteq \mathbb{N}^d$ such that $\gcd(a_1,a_2)=1$, $d_j\in \mathbb{N}\setminus \{0\}$ for all $j\in \{1,\ldots,d\}\setminus \{i\}$, $\mathbf{x}\notin \langle A\setminus \{\mathbf{x}\}\rangle$ and $\mathbf{x}\notin \{\lambda \mathbf{e}_i, \mathbf{e}_i+\beta_j\mathbf{e}_j \mid \lambda,\beta \in \mathbb{N}\}$. Let $S=\langle A \rangle $, then the following holds:
    \begin{enumerate}
        \item If $d\geq 3$ then $|\operatorname{FA}(S)|\geq 2$.
        \item If $d=2$ and $a\neq3$ then $|\operatorname{FA}(S)|\geq 2$.
        \item If $d=2$ and $a_1=3$, then $|\operatorname{FA}(S)|=1$ if and only if $\mathbf{x}=(2,d_j)$ with $j\in \{1,2\}\setminus \{i\}$. %Moreover, when $|\operatorname{FA}(S)|=1$, we get $\operatorname{FA}(S)=\{(\operatorname{F}(\langle a_1,a_2\rangle), d_2-1)\}$.
        % \item If $d=2$ and $a_1=3$, then \begin{enumerate}
        % \item  if $\mathbf{x}\neq (2,d_2)$, then $|\operatorname{FA}(S)|\geq 2$;
        % \item if $\mathbf{x}=(2,d_2)$, then $\operatorname{FA}(S)=\{(\operatorname{F}(\langle a_1,a_2\rangle), d_2-1)\}$.
        %\end{enumerate}
    \end{enumerate}
\end{lemma}

\begin{proof}
    Without loss of generality, we can assume that $i=1$. Let us denote $T=\langle a_1,a_2\rangle$ and $S'=\langle A\setminus \{\mathbf{x}\}\rangle$. Observe that if $\mathbf{s}\in S'$, then $\mathbf{s}\in S$.

    \noindent Suppose first $a_1=2$. By the proof of \cite[Lemma 4.1]
    {cisto2020generalization} we have that 
    $$\cH(S')=\left\lbrace h\mathbf{e}_1+\sum_{j\neq 1} \lambda_j\mathbf{e}_j\mid h\in \cH(T), \lambda_j\in \{0,1,\ldots,d_j-1\}, j\in \{2,\ldots,d\}\right \rbrace.$$
Therefore, $\mathbf{x}=x_1\mathbf{e}_1+\sum_{j\neq i}x_j\mathbf{e}_j$ with $x_1\notin T$ and $x_j<d_j$ for all $j\in \{2,\ldots,d\}$. 
    %Also, by our assumptions on $\mathbf{x}$, there exists $k\neq 1$ such that $x_k>0$. In particular, $\mathbf{x}-\mathbf{e}_k\in \cH(S')$. Moreover $\mathbf{x}-\mathbf{e}_k\in \cH(S)$, since if it belongs to $S$, then $\mathbf{x}\in S'$ because $\mathbf{x}$ is not involved in any factorization of $\mathbf{x}-\mathbf{e}_k$ in $S$. 
    Consider $\mathbf{f}=\operatorname{F}(T)\mathbf{e}_1+(x_k-1)\mathbf{e}_k+\sum_{j\neq 1,k}(d_j-1)\mathbf{e}_j$. We have that $\mathbf{f}\in \cH(S)$. In fact, since $\mathbf{f}\in \cH(S')$, the only possibility to $\mathbf{f}\in S$ is $\mathbf{f}=\mathbf{x}+\mathbf{s}$ for some $\mathbf{s}\in S$, but this is a contradiction since the $k$-th coordinate of $\mathbf{x}+\mathbf{s}$ is bigger than $x_k-1$. Now we prove that $\mathbf{f}\in \operatorname{FA}(S)$. Consider $\mathbf{y}=\sum_{j=1}^d y_j\mathbf{e}_j\in \mathbb{N}^d\setminus \{\mathbf{0}\}$. It is clear that if $y_1>0$ or $y_j>0$ for some $j\neq k$, then $\mathbf{f}+\mathbf{y}\in S'\subseteq S$. So, it is sufficient to show that $\mathbf{f}+y_k\mathbf{e}_k\in S$ for every $y_k>0$. Now,  
    $$\mathbf{f}+y_k\mathbf{e}_k= \mathbf{x}+(\operatorname{F}(T)-x_1)\mathbf{e}_1+(y_k-1)\mathbf{e}_k+\sum_{j\neq 1,k}(d_j-1-x_j)\mathbf{e}_j \in S.$$
    Since $T$ is symmetric and $x_1\notin T$, we obtain $\operatorname{F}(T)-x_1 \in T$. Hence $\mathbf{f}\in \operatorname{FA}(S)$. Now consider $\mathbf{g}=\mathbf{e}_1+x_k\mathbf{e}_k$. Observe that $\mathbf{g}\notin S'$ and if $\mathbf{g}\in S$, then $\mathbf{g}=\mathbf{x}+\mathbf{s}$ for some $\mathbf{s}\in S$. This forces that $\mathbf{x}=\mathbf{e}_1+x_k\mathbf{e}_k$, this is a contradiction to our hypothesis. Therefore, $\mathbf{g}\in \cH(S)$ and it is not comparable to $\mathbf{f}$ with respect to $\leq_{\mathbb{N}^d}$. Thus, we obtain $|\operatorname{FA}(S)|\geq 2$.

    \noindent Suppose $a_1\geq 3$. By the proof of \cite[Lemma 4.1]
    {cisto2020generalization} we have that $\mathbf{f}=\operatorname{F}(T)\mathbf{e}_1+\sum_{j\neq 1}(d_j-1)\mathbf{e}_j$ is a maximal gap of $S'$ with respect to $\leq_{\mathbb{N}^d}$. Let us denote $\mathbf{x}=x_1\mathbf{e}_1+\sum_{j\neq i}x_j\mathbf{e}_j$. In particular, since $\mathbf{x}\notin S'$, we have $x_1\notin T$. Suppose there exists $k\in \{2,\ldots,d\}$ such that $x_k>d_k$. In this case $\mathbf{f}\notin S$, otherwise we have $\mathbf{f}=\mathbf{x}+\mathbf{s}$ for some $\mathbf{s}\in S$ and this not possible because the $k$-th coordinate of $\mathbf{f}$ is smaller than $x_k$. Therefore, $\mathbf{f}\in \cH(S)$ and it is also a maximal gap of $S$ with respect to $\leq_{\mathbb{N}^d}$ since it is a maximal gap for $S'$. Moreover $\mathbf{x}-\mathbf{e}_k\in \cH(S)$, since if it belongs to $S$, then $\mathbf{x}\in S'$ because $\mathbf{x}$ is not involved in any factorization of $\mathbf{x}-\mathbf{e}_k$ in $S$. In particular, we obtain $\mathbf{x}-\mathbf{e}_k\in \cH(S)$ and it is not comparable to $\mathbf{f}$ with respect to $\leq_{\mathbb{N}^d}$. Hence, $|\operatorname{FA}(S)|\geq 2$. Now, suppose that $x_j\leq d_j$ for all $j\in \{2,\ldots,d\}$. Observing that $\operatorname{F}(T)\mathbf{e}_1\in \cH(S)$, we have that there exists a maximal gap of $S$ of kind $\mathbf{f}=\operatorname{F}(T)\mathbf{e}_1+\sum_{j\neq 1}\alpha_j\mathbf{e}_j$. We show that $\alpha_j<d_j$ for all $j\in \{2,\ldots,d\}$. If $\alpha_k\geq d_k$ for some $k\in \{2,\ldots,d\}$, we get 
    $$\mathbf{f}=(\mathbf{e}_1+d_k\mathbf{e}_k)+(\operatorname{F}(T)-1)\mathbf{e}_1+\sum_{j\neq 1,k}\alpha_j\mathbf{e}_j+(\alpha_k-d_k)\mathbf{e}_k\in S,$$ since $\operatorname{F}(T)-1\in T$, a contradiction. Consider the elements $\mathbf{y}_k=2\mathbf{e}_1+d_k\mathbf{e}_k$ for $k\in \{2,\ldots,d\}$. Observe that $\mathbf{y}_k\in \cH(S')$ for all $k\in \{2,\ldots,d\}$. Suppose there exists $\ell\in \{2,\ldots,d\}$ such that $\mathbf{y}_\ell \in S$. Hence, we have that $\mathbf{y}_\ell=\mathbf{x}+\mathbf{s}$ for some $\mathbf{s}\in S$. This means that $\mathbf{x}=2\mathbf{e}_1+x_\ell\mathbf{e}_\ell$ for some $0<x_\ell\leq d_\ell$. In particular, if $d>2$, then there exists $\ell'\in \{2,\ldots,d\}\setminus \{\ell\}$ such that $\mathbf{y}_{\ell'}\in \cH(S)$ and $\mathbf{y}_{\ell'}$ is not comparable to $\mathbf{f}$. Hence, if $d>2$ we get $|\operatorname{FA}(S)|\geq 2$. So, we now assume $d=2$. For simplicity, we denote $S=\langle (a_1,0),(a_2,0),(0,1),(1,d_2),(2,x_2)\rangle$ with $0<x_2\leq d_2$ and $a_1\geq 3$. If $a_1>3$ then $(3,x_2-1+d_2)\in \cH(S)$ and it is not comparable with $\mathbf{f}=(\operatorname{F}(T),\alpha_2)$, that is, $|\operatorname{FA}(S)|\geq 2$. It remains to consider $a_1=3$. Note that, in this case, $\alpha_2<x_2$, otherwise we obtain $(\operatorname{F}(T),\alpha_2)=(\operatorname{F}(T)-2,0)+(2,x_2)+(0,\alpha_2-x_2)\in S$. In particular, if $x_2<d_2$ then $(1,d_2-1)\in \cH(S)$ and it is not comparable to $\mathbf{f}=(\operatorname{F}(T),\alpha_2)$. Now, consider $x_2=d_2$. In this case it is not difficult to see that $\cH(S)=\{(h,t)\mid h\in \cH(T), 0\leq t\leq d_2-1\}$. In particular, $\operatorname{FA}(S)=\{(\operatorname{F}(T),d_2-1)\}$. This completes the proof. 
    \end{proof}

    \section{Main result and some questions}

With the help of the above lemmas, one can classify the symmetric generalized numerical semigroups of embedding dimension $2d+1$ in $\mathbb{N}^d.$  If $d=1$, then $S$ is a numerical semigroup of embedding dimension $3$. Let $S = \langle a_1,a_2,a_3 \rangle$ be a numerical semigroup of embedding dimension three. Define $c_i = \mathrm{min}\{\lambda_i \in \mathbb{N} \mid \lambda_i \neq 0, \lambda_ia_i \in \langle \{a_1,a_2,a_3\} \setminus \{a_i\} \rangle \}.$ From \cite{Herzog}, we know that $S = \langle a_1,a_2,a_3 \rangle$ is symmetric if and only if $c_ia_i = c_ja_j$ for some $i \neq j.$ Therefore, we consider $d>1$, and have the following:

\begin{theorem}\label{thm:symmteric2d+1}
Let $d>1$, and $S\subseteq \mathbb{N}^d$ be a generalized numerical semigroup with embedding dimension $2d+1$. Then $S$ is symmetric if and only if $|\operatorname{FA}(S)|=1$.  Moreover, in this case, sets of symmetric and almost symmetric generalized numerical semigroups coincide.
% Let $S$ be a GNS as in Lemma~\ref{lem:class1} or Lemma~\ref{lem:class3}. Then $S$ is symmetric if and only if $|\operatorname{FA}(S)|=1$.    
\end{theorem}
    \begin{proof}
    It is always true that if $S$ is symmetric, then $|\operatorname{FA}(S)|=1$. Now assume $|\operatorname{FA}(S)|=1$. So, $S$ is as in Lemma~\ref{lem:class1} or Lemma~\ref{lem:class3}. If $S$ is as in Lemma~\ref{lem:class1} and $|\operatorname{FA}(S)|=1$, then we have $S\subseteq \mathbb{N}^2$ and $S=\langle (3,0),(a_2,0),(a_3,0),(0,1),(1,d_2)\rangle$ with $a_2\equiv 1 \mod 3$, $a_3=a_2+1$ and $d_2\in \mathbb{N}\setminus \{0\}$. Let us denote by $T$ the numerical semigroup generated by the set $\{3,a_2,a_3\}$. In this case, considering also the arguments in the proof of Lemma~\ref{lem:class1}, we have that $\operatorname{FA}(S)=\{(\operatorname{F}(T), 2d_2-1)\}$. Moreover, the set of gaps is $\cH(S)=\{(h,j)\mid h\in \cH(T), h\equiv 1 \mod 3, 0\leq j\leq d_2-1\}\cup \{(h,j)\mid h\in \cH(T), h\equiv 2 \mod 3, 0\leq j\leq 2d_2-1\} $. Observe that $\operatorname{F}(T)=a_2-2$ and $\cH(T)=\{h<a_2-1 \mid h\text{ is not divisible by }3 \}$. So, $|\cH(T)|= \frac{2}{3}(a_2-1)$. Observe also that the sets $\{h\in \cH(T)\mid h\equiv 1 ~\text{mod}~ 3\}$ and $\{h\in \cH(T)\mid h\equiv 2 ~\text{mod}~ 3\}$ have the same cardinality. Therefore, we get $$2|\cH(S)|= \frac{2}{3}(a_2-1)\cdot d_2+\frac{2}{3}(a_2-1)\cdot 2d_2=2d_2(a_2-1).$$
    Considering that $\operatorname{FA}(S)=\{(a_2-2,2d_2-1)\}$, we obtain that $S$ is symmetric by Theorem~\ref{symcrit}. If $S$ is as in Lemma~\ref{lem:class3} and $|\operatorname{FA}(S)|=1$, then we have $S\subseteq \mathbb{N}^2$ and $S=\langle (3,0),(a_2,0),(0,1),(1,d_2), (2,d_2)\rangle$. Let us denote by $T$ the numerical semigroup generated by the set $\{3,a_2\}$. Considering the arguments in the proof of Lemma~\ref{lem:class3}, we have $\cH(S)=\{(h,t)\mid h\in \cH(T), 0\leq t\leq d_2-1\}$ and $\operatorname{FA}(S)=\{(\operatorname{F}(T),d_2-1)\}$. Since $T$ is a numerical semigroup of embedding dimension 2, we know that $|\cH(T)|=\frac{\operatorname{F}(T)+1}{2}$ (see \cite{sylvester1884mathematical} or \cite[Proposition 2.13]{rosales2009numerical}). Therefore, $2|\cH(S)|=(\operatorname{F}(T)+1)d_2$ and we obtain that $S$ is symmetric by Theorem~\ref{symcrit}.

    Moreover, recall that if $S$ is almost symmetric, then $S$ is symmetric. Conversely, if $S$ is almost symmetric then $|\operatorname{FA}(S)|=1$. This concludes the proof.
    \end{proof}

% Now, one can also classify the almost symmetric generalized numerical semigroups of embedding dimension $2d+1$ in $\mathbb{N}^d$. If $d=1$, then $S$ is a numerical semigroup of embedding dimension $3$. In this case, almost symmetric semigroups are either pseudo-symmetric or symmetric, as the type of a numerical semigroup generated by three elements is at most two. For a complete classification of symmetric numerical semigroups of embedding dimension three, one can see \cite{Herzog}. Also, for a complete classification of pseudo-symmetric numerical semigroups of embedding dimension three, one can see \cite[Theorem 15]{pedroEmbedding3}. Therefore, we consider $d > 1$, and have the following:

% \begin{corollary}
%     Let $d>1$, and $S\subseteq \mathbb{N}^d$ be a generalized numerical semigroup of embedding dimension $2d+1$. Then, $S$ is almost symmetric if and only if $S$ is symmetric.
% \end{corollary}
% \begin{proof}
%     A symmetric generalized numerical semigroup is trivially almost symmetric, while if $S$ is an almost symmetric generalized numerical semigroup, then $|\operatorname{FA}(S)|=1$. On the other hand, in this case, $S$ is symmetric by Theorem \ref{thm:symmteric2d+1}. 
% \end{proof}

\begin{remark}
    {\rm When $d >1$, one gets that there does not exist any generalized numerical semigroup of embedding dimension $2d$ or $2d+1$, which is almost symmetric but not symmetric. But, when the embedding dimension is equal to $2d+2$, there exist almost symmetric generalized numerical semigroups that are not symmetric, as shown in the following example.}
\end{remark}

\begin{example} \rm
    Let $S=\langle (0, 3 ), (0, 4 ), ( 0, 5 ), (1, 0 ), ( 4, 1 ), ( 7, 2 )\rangle \subseteq \mathbb{N}^2$. Then $S$ is a generalized numerical semigroup of embedding dimension $6$, having set of gaps $\mathcal{H}(S)=\{ ( 0, 1 ), (0, 2), (1, 1), (1, 2 ), (2, 1 ), (2, 2 ), (3, 1 ), (3, 2 ), (4, 2 ), (5, 2 ), (6, 2 )\}$. It is not difficult to compute $\operatorname{PF}(S)=\{(3,1),(6,2)\}$. In particular, $S$ is almost symmetric but not symmetric.
\end{example}

We have seen that the type of almost symmetric generalized numerical semigroups with embedding dimension $2d+1$ is bounded by $2$. In general, it is known that the type of generalized numerical semigroups of fixed embedding dimension $e$ may not be a bounded function of $e$ (see \cite[Section 6]{op4}). But the same question is unknown for almost symmetric generalized numerical semigroups. In particular, fix an integer $e$ and consider the set of all almost symmetric generalized numerical semigroups with embedding dimension $e$. Does an upper bound exist for the type of these semigroups in terms of $e$?
%for all $d \geq 1$ and all values of embedding dimension $e \geq 2d$, we have already given in \cite{preprint} \textcolor{red}{We could add the reference of part1} a class of examples conveying that the type of a generalized numerical semigroup in $\mathbb{N}^d$ may not be a bounded function of embedding dimension $e.$ 
%The following question is open-ended for almost symmetric numerical semigroups, and we ask it in more generality.
%\begin{question}
    %Fix an integer $e \geq 6$ and consider the set of all almost symmetric generalized numerical semigroups with embedding
%dimension $e$. Does an upper bound exist for the type of these semigroups?
%\end{question}
The answer to this question is affirmative for small embedding dimensions, in particular for $e =2d, 2d+1$, and for $d =1$, the answer is affirmative for $e \leq 5$ (see \cite{moscarielloedim4}, \cite{moscarielloedim5}). We computationally examined the above question for generalized numerical semigroups with embedding dimension $e=2d+2$. In particular, with the help of the \texttt{GAP}\cite{GAP} package \texttt{numericalsgps}\cite{NumericalSgps}, we explored the set of almost symmetric generalized numerical semigroups in $\mathbb{N}^2$ such that the cardinality of the set of gaps is smaller or equal than $14$ and the set of almost symmetric generalized numerical semigroups in $\mathbb{N}^3$ such that the cardinality of the set of gaps is smaller or equal than $10$. In this case, we found that the type is bounded by $2$. 

\noindent  More specifically, we identified the generalized numerical semigroups up to the permutation of coordinates for the computations. It is not difficult to see that if $S$ is a generalized numerical semigroup and $S'$ is a generalized numerical semigroup obtained from $S$ by a permutation of coordinates, then $S$ is almost symmetric (symmetric, pseudo-symmetric) if and only if $S'$ is almost symmetric (symmetric, pseudo-symmetric), and the two semigroups have the same type. A procedure to obtain the set of generalized numerical semigroups up to permutation of coordinates and the fixed genus is provided in \cite{cisto2024generalized}. After this computation, it is not difficult to obtain all the generalized numerical semigroups considering all possible permutations of coordinates. This trick allowed us to reach bigger values of the genus in a reasonable time (especially in dimension 3). Implementing this procedure in \texttt{GAP}, we obtained the computational results described in Table~\ref{table1}, where the following notations are used:  

\begin{itemize}
\item $\mathcal{AS}_{g,d}(e)$ is the set of almost symmetric generalized numerical semigroups in $\mathbb{N}^d$ such that the cardinality of the set of gaps is $g$ and the embedding dimension is equal to $e$.
\item $\mathcal{S}ym_{g,d}(e)$ is the set of symmetric generalized numerical semigroups in $\mathbb{N}^d$ such that the cardinality of the set of gaps is $g$ and the embedding dimension is equal to $e$.
\item $\mathcal{PS}ym_{g,d}(e)$ is the set of pseudo-symmetric generalized numerical semigroups in $\mathbb{N}^d$ such that the cardinality of the set of gaps is $g$ and the embedding dimension is equal to $e$.
\end{itemize}

\begin{table}[h!]
\begin{tabular}{|c|cccccc|}
\toprule
$g$  & $|\mathcal{AS}_{g,2}(6)|$  & $|\mathcal{S}ym_{g,2}(6)|$ & $|\mathcal{PS}ym_{g,2}(6)|$ & $|\mathcal{AS}_{g,3}(8)|$  & $|\mathcal{S}ym_{g,3}(8)|$ & $|\mathcal{PS}ym_{g,3}(8)|$\\
\midrule

  1& 0& 0& 0& 0& 0& 0\\
  2& 2& 0& 2& 0& 0& 0\\
  3& 4& 2& 2& 9& 9& 0\\
  4& 0& 0& 0& 3& 3& 0\\
  5& 8& 4& 4& 0& 0& 0\\
  6& 8& 6& 2& 27& 27& 0\\
  7& 4& 4& 0& 3& 3& 0\\
  8& 10& 4& 6& 6& 6& 0\\
  9& 10& 8& 2& 27& 27& 0\\
  10& 4& 4& 0& 3& 3& 0\\
  11& 10& 6& 4& 0& 0 & 0\\
  12& 12& 10& 2& & & \\
  13& 6& 4& 2& & & \\
  14& 8& 4& 4& & &\\
\bottomrule
\end{tabular}
\caption{Computational data}
\label{table1}
\end{table} 

\begin{example} \rm
Using the procedure described in \cite{cisto2024generalized}, it is possible to compute all generalized numerical semigroups in $\mathbb{N}^2$, up to permutation of coordinates and having genus equal to $3$. The total number of these semigroups is $12$. Among them, there are $2$ almost symmetric with embedding dimension equal to $2d+2=6$, which are the following: 
\begin{itemize}
\item $S_1=\mathbb{N}^2\setminus \{(0,1),(0,2),(0,4)\}= \langle (1,0), (0,3), (1,1), (1,2), (0,5), (0,7) \rangle$.
\item $S_2=\mathbb{N}^2\setminus \{(0,1),(1,0),(1,2)\}=\langle (0,3), (1,1), (0,2), (3,0), (2,0), (2,1)\rangle$.
\end{itemize}
It is not difficult to check that $S_1$ is pseudo-symmetric and $S_2$ is symmetric. Considering all permutation of coordinates, the following semigroups should be added:
\begin{itemize}
\item $S'_1=\mathbb{N}^2\setminus \{(1,0),(2,0),(4,0)\}=\langle (0,1), (3,0), (1,1), (2,1), (5,0), (7,0) \rangle $.
\item $S'_2=\mathbb{N}^2\setminus \{(0,1),(1,0),(2,1)\}=\langle (0,3), (1,1), (0,2), (3,0), (2,0), (1,2)\rangle$
\end{itemize}
Therefore, we have $\mathcal{AS}_{3,2}(6)=\{S_1,S'_1,S_2,S'_2\}$, $\mathcal{S}ym_{3,2}(6)=\{S_2,S'_2\}$ and $\mathcal{PS}ym_{3,2}(6)=\{S_1,S'_1\}$.

\end{example}

The data in Table~\ref{table1} show that all almost symmetric are symmetric or pseudo-symmetric, with the parameters considered. This means that the type is bounded by 2. Having the positive answer for $d =1$ (i.e., for numerical semigroups) and considering our computational experiments, we pose the following conjecture:

\begin{conjecture} \label{conj:2d+2}
    Let $S \subseteq \mathbb{N}^d$ be an almost symmetric generalized numerical semigroup of embedding dimension $2d+2$. Then, the type of $S$ can be at most three.
\end{conjecture}

{\it Acknowledgment.} This work was done during the first author's visit to the University of Messina, Italy, in May 2024. The first author would like to extend sincere thanks to the University of Messina and the group GNSAGA of Istituto Nazionale di Alta Matematica (INdAM), Italy, for their support.

% If using bibtex:

%\nocite{*}
\bibliographystyle{plain}
\bibliography{biblio}

%\begin{thebibliography}{99}
%		
%		\addcontentsline{toc}{chapter}{\bibname}
%		
%		\bibitem{C} Citation .
%		
%				
%\end{thebibliography}
	\end{document}